\documentclass[12pt]{article}
\usepackage{amssymb,amsmath,latexsym}

\newtheorem{theorem}{Theorem}[section]

\newtheorem{remark}[theorem]{Remark}

\newenvironment{proof}{\noindent{\sc Proof.}}{\hfill\qed}
\newcommand{\qed}{\quad\lower0.05cm\hbox{$\Box$}}

\newcommand{\D}{\mathcal{D}}
\newcommand{\N}{\mathcal{N}}
\newcommand{\Na}{\mathbb{N}}
\newcommand{\R}{{\mathbb R}}

\newcommand{\U}{{\mathcal U}}
\newcommand{\B}{{\mathcal B}}

\newcommand{\abs}[1]{\mid \hspace{-.1cm} #1 \hspace{-.1cm}\mid}

\title{Translation-invariant generalized topologies induced by probabilistic norms}
\author{
Bernardo Lafuerza-Guill\'{e}n\\ {\small Departamento de
Estad\'{\i}stica y Matem\'{a}tica Aplicada, Universidad de
Almer\'\i a,}\\ {\small E--04120 Almer\'\i a, Spain, e-mail: {\tt
blafuerz@ual.es}}\\
\and Jos\'{e} L. Rodr\'{\i}guez \\ {\small \'{A}rea de
Geometr\'{\i}a y Topolog\'{\i}a, Universidad de Almer\'\i a,}\\
{\small E--04120 Almer\'\i a, Spain, e-mail: {\tt
jlrodri@ual.es}}}
\date{\empty}

\begin{document}
\maketitle
\begin{abstract}
In this paper we consider probabilistic normed spaces as defined
by Alsina, Sklar, and Schweizer, but equipped with non necessarily
continuous triangle functions. Such spaces endow a generalized
topology that is Fr\'echet-separable, translation-invariant and
countably generated by radial and circled $0$-neighborhoods.
Conversely, we show that such generalized topologies are
probabilistically normable.


\end{abstract}

\section{Introduction}
Probabilistic normed spaces (briefly, PN spaces) were first
defined by \v Serst\-nev in the early sixties (see \cite{Serst}),
deriving to a fruitful theory concordant with the theory of
ordinary normed spaces. Twenty years later, Alsina, Schweizer, and
Sklar gave in \cite{alsina1} a quite general definition of PN
space, based on the definition of Menger's betweenness in
probabilistic metric spaces; see \cite[p.~232]{SS}.

We here consider PN spaces in which the involved triangle
functions are non necessarily continuous. The known underlying
strong topology (exploited in \cite{alsina1}) is no longer a
genuine topology but a generalized topology in the sense of
Fr\'{e}chet. The analogous situation for probabilistic metric
spaces was treated by H\"{o}hle in~\cite{Hoh}, where he showed
that all generalized topologies which are Fr\'{e}chet-separated
and first-numberable are induced by certain probabilistic metrics.
The main result of this paper is a similar result for
probabilistic norms:

\medskip
\noindent {\bf Theorem 4.2} {\it Let $T$ be a $t$-norm such that
$\sup_{0\leq x<1} T(x,x)<1$. Suppose that $T(x,y)\leq xy$,
whenever $x, y<\delta$, for some $\delta>0$. A
Fr\'{e}chet-separated, translation-invariant, generalized topology
$(\U_p)_{p\in S}$ on a real vector space $S$ is derivable from a
Menger PN space $(S,\nu,\tau_T,\tau_{T^*})$, if and only
$\U_{\theta}$ admits a countable base of radial and circled
subsets, where $\theta$ is the origin of $S$.}

\bigskip

We think that a similar result could be interesting for fuzzy
normed spaces in the sense of Felbin \cite{Fel}, but allowing
non-continuity of the $t$-norms, and $t$-conorms involved in the
fuzzy structure.

In \cite{LR} the authors use this generalized topology to define
bounded subsets in PN spaces (with non necessarily continuous
triangle functions) and study its relationship with $\D$-bounded
subsets (a concept which is defined in probabilistic terms).




\section{PM and PN spaces}

Recall from \cite{SS} and \cite{alsina1} some definitions on
probabilistic metric and probabilistic normed spaces.

As usual $\Delta^+$ denotes the set of distance distribution
functions (briefly, a d.d.f.), i.e.~distribution functions with
$F(0)=0$, endowed with the metric topology given by the modified
Levy-Sybley metric $d_L$ (see 4.2 in \cite{SS}).
Given a real number $a$, $\varepsilon_a$ denotes the distribution
function defined as $\varepsilon_a(x)=0$ if $x\leq a$ and
$\varepsilon_a(x)=1$ if $x> a$. Hence, $\R^+$ can be viewed as a
subspace of $\Delta^+$. A triangle function $\tau$ is a map from
$\Delta^+\times \Delta^+ \to \Delta^+$ which is commutative,
associative, nondecreasing in each variable and has
$\varepsilon_0$ as the identity. Such functions give rise to all
possible extensions of the sum of real numbers, so that (M3) below
corresponds to the triangle inequality.


A {\it probabilistic metric space} (briefly, a PM space) is a
triple $(S,F,\tau)$ where $S$ is a non-empty set, $F$ is a map
from $S\times S \to \Delta^+$, called the probabilistic metric,
and $\tau$ is a triangle function, such that:
\begin{itemize}
\item[(M1)] $F_{p,q}=\varepsilon_0$ if and only if $p=q$.
\item[(M2)] $F_{p,q}=F_{q,p}$. \item[(M3)]
$F_{p,q}\geq\tau(F_{p,r},F_{r,q})$.
\end{itemize}
When only (M1) and (M2) are required, it the pair $(S,F)$ is said
to be a {\it probabilistic semi-metric space} (briefly, PSM
space).

A {\it PN space}
 is a
quadruple $(S,\nu,\tau,\tau^*)$ in which $S$ is a vector space
over $\R$, the {\it probabilistic norm} $\nu$ is a map $S\to
\Delta^+$, $\tau$ and $\tau^*$ are
triangle functions\footnote{In the definition of PN space given in
\cite{alsina1} the triangle functions are assumed to be
continuous} such that the following conditions are satisfied for
all $p$, $q$ in $S$:
\begin{itemize}
\item[(N1)] $\nu_p=\varepsilon_0$ if and only if $p=\theta$, where
$\theta$ is the origin of $S$. \item[(N2)] $\nu_{-p}=\nu_p$.
\item[(N3)] $\nu_{p+q}\geq\tau(\nu_p,\nu_q)$. \item[(N4)]\,
$\nu_p\leq\tau^*(\nu_{\lambda p},\nu_{(1-\lambda)p})$ for every
$\lambda\in [0,1]$.
\end{itemize}

Observe that every PN space $(S,\nu,\tau,\tau^*)$ is a PM space,
where $F_{p,q}:=\nu_{p-q}$.

 Recall that a {\it $t$-norm} is a binary operation on $[0,1]$ that is commutative, associative,
nondecreasing in each variable, and has $1$ as identity. Dually, a
$t$-conorm is a binary operation on $[0,1]$ that is commutative,
associative, non-increasing in each variable, and has $0$ as
identity. If $T$ is a $t$-norm, its associated $t$-conorm $T^*$ is
defined by $T^*(x,y):= 1-T(1-x, 1-y)$. Given a $t$-norm $T$ one
defines the functions $\tau_T$ and $\tau_{T^*}$ by
$$\tau_T(F,G)(x):=\sup\{T(F(s), G(t)): s+t=x\},$$ and
$$\tau_{T^*}(F,G)(x):=\inf \{T^*(F(s), G(t)): s+t=x\}.$$ Recall that if $T$ is
left-continuous then $\tau_T$ is a triangle function
\cite[p.~100]{SS}, although this is not necessary;
 For example, if $Z$ denotes the minimum $t$-norm, defined as
$Z(x,1)=Z(1,x)=x$ and $Z(x,y)=0$ elsewhere, then $\tau_Z$ is a
triangle function which is not continuous.

 A  {\it \v Serstnev PN space} is a PN space
$(V,\nu,\tau,\tau^*)$ where $\nu$ satisfies the following \v
Serstnev condition:
\begin{itemize}
\item[(\v {S})] \quad  $\nu_{\lambda p}(x)=
\nu_p\left(\frac{x}{\mid \lambda\mid }\right)$, for all $x\in
\R^+$, $p\in V$ and $\lambda \in \R\setminus \{0\}$.
\end{itemize}
It turns out that (\v S) is equivalent to have (N2) and
\begin{equation}
\label{characterization-serstnev} \nu_p=\tau_M(\nu_{\lambda
p},\nu_{(1-\lambda)p}),
\end{equation}
for all $p\in V$ and $\lambda\in [0,1]$ (see
\cite[Theorem~1]{alsina1}), where $M$ is the $t$-norm defined as
$M(x,y)=\min\{x,y\}$.

Let $T$ be a $t$-norm. A {\it Menger PM space under $T$} is a PM
space of the form $(S,F,\tau_T)$. Analogously, a {\it Menger PN
space}
is a PN space
of the form $(S, \nu, \tau_T, \tau_{T^*})$. Note that every metric space $(S,d)$
is a Menger space $(S,F,\tau_M)$ where $F_{p,q}=\varepsilon_{d(p,q)}$.
Analogously, every normed space $(S,{\parallel \cdot \parallel})$ is a Menger and
\v Serstnev PN space $(S,\nu,\tau_M, \tau_{M})$
where $\nu_p=\varepsilon_{\parallel p \parallel}$.

\section{Probabilistic metrization of generalized topologies} \label{Hohle-paper}
In \cite{Hoh} H\"ohle solved a problem posed by Thorp about the probabilistic metrization
of generalized topologies. We recall some definitions and results that we shall used
in the next section.

Let $S$ be a (non-empty) set. A {\it generalized topology (of type
$V_D$)} on $S$ is a family of subsets $(\U_p)_{p\in S}$, where
$\U_p$ is a filter on $S$ such that $p\in U$ for all $U\in \U_p$
(see e.g. \cite[p.~38]{SS}, \cite[p.~22]{AK}). Elements of $\U_p$
are called {\it neighborhoods} at $p$. Such a generalized topology
is called {\it Fr\'{e}chet-separable\,} if $ \bigcap_{U\in
\U_p}U=\{p\}$.

A {\it generalized uniformity} $\U$ on $S$ is a filter on $S\times
S$ such that every $V\in \U$ contains the diagonal $\{(p,p): p\in
S\}$, and for all $V\in \U$, we have that $V^{-1}:= \{(q,p):
(p,q)\in V\}$ also belongs to $\U$. Elements of $\U$ are called
{\it vicinities}. Every generalized uniformity $\U$ induces a
generalized topology as follows: for $p\in S$,
\begin{equation}
\label{GU} \U_p:= \{U\subseteq S \mid \exists V\in \U: U\supseteq
\{q\in S \mid (p,q)\in V\}\}.
\end{equation}
A uniformity $\U$ is called {\it Hausdorff-separated} if the
intersection of all vicinities is the diagonal on $S$.
Theorem 1 in \cite{Hoh} claims:

\begin{theorem}[H\"ohle] Every Fr\'{e}chet-separable generalized
topology $(\U_p)_{p\in S}$ on a given set $S$ is derivable from a
Hausdorff-separable generalized uniformity $\U$ in the sense of
(\ref{GU}).\qed
\end{theorem}

 Let $(S,F)$ be a PSM space.
Consider the system $(\N_p)_{p\in S}$, where $\N_p=\{N_p(t): t>0\}$ and
$$N_p(t):=\{q\in S: F_{p,q}(t)> 1-t\}.$$
This is called the strong neighborhood system. If we define
$\delta(p,q):= d_L(F_{p,q},\varepsilon_0)$, then $\delta$ is a
semi-metric on $S$ (i.e.~it may not satisfy the triangle
inequality of the standard metric axioms), and $N_p(t)=\{q:
d_L(F_{p,q}, \varepsilon_0)<t\}$. Clearly $p\in N$ for every $N\in
\N_p$, and the intersection of two strong neighborhoods at $p$ is
a strong neighborhood at $p$. Furthermore, $\N_p$ admits a
countable filter base given by  $\{N_p(1/n): n\in \Na\}$, hence
the strong neighborhood system is first-numerable. The above
explanation yields the following fact (see more details in
\cite[p. 191]{SS}):

\begin{theorem} Let $(S,F)$ be a PSM space, then the strong
neighborhood system defines a generalized topology of type
$V_D$ which is Fr\'echet-separable and first-numerable. \qed
\end{theorem}
This generalized topology is called the {\it strong generalized
topology} of the PSM space $(S,F)$.

The main result in \cite{Hoh} is the following.

\begin{theorem}[H\"ohle]
\label{Hohle main theorem} Let $T$ be a $t$-norm such
that $\sup_{0\leq x<1} T(x,x)<1$. A Fr\'echet-separated
 generalized topology $(\U_p)_{p\in S}$ on a
set $S$ is derivable from a Menger PM space $(S, F, T)$ if
and only if there exists a Hausdorff-separated, generalized uniform structure $\U$
having a countable filter base, such that $\U$ is compatible with $(\U_p)_{p\in S}$. \qed
\end{theorem}

\begin{remark} {\rm If $(S,F,\tau)$ is a PM space with $\tau$ continuous, then the
associated generalized topology is in fact a topology (i.e.~it satisfies that
for every $N\in \N_p$ and $q\in N$, there is a $N'\in \N_q$ such
that $N'\subseteq N$). This topology is called {\it the strong
topology}. Because of (M1)   this topology is Hausdorff. Since it
is first-numerable and uniformable, one has that it is metrizable
(see \cite[Theorem 12.1.6]{SS}).

Conversely, if $\sup_{0\leq x<1} T(x,x)=1$, then a
Fr\'echet-separated, uniformable, topology is derivable from a
Menger space $(S,F,T)$ if and only if there exists a Hausdorff
uniformity $U$ on $S$ having a countable filter base
(\cite{Hoh}).}
\end{remark}
\section{Translation-invariant generalized topologies}

Assume now that $S$ is a vector space over $\R$. A generalized
topology $(\U_p)_{p\in S}$ on $S$ is said to be  {\it
translation-invariant} if for all $U\in \U_p$ and $q\in S$, we
have $q+U\in \U_{p+q}$. Consequently, a translation-invariant
generalized topology is uniquely determined by the neighborhood
system $\U_\theta$ at the origin $\theta$ of $S$. In this case,
the generalized uniformity from which one can derive the
generalized topology is:
$$
\U:= \{V\subseteq S\times S \mid \exists U\in \U_\theta :
V\supseteq \{(p,q)\mid p-q\in U\}\}.
$$
Recall that a subset $U$ of a vector space is called {\it radial} if $-U=U$; it is
called {\it circled (or balanced)} if  $\lambda U\subset U$ for all $\abs{\lambda}\leq
1$.

\begin{theorem}
\label{maintheorem} Every PN space $(S, \nu, \tau, \tau^*)$ admits
a generalized topology $(\U_p)_{p\in S}$ of type $V_D$ which is
Fr\'{e}chet-separable, translation-invariant, and
countably-generated by radial and circled $\theta$-neighborhoods.
\end{theorem}
\begin{proof}
Let $(S,\nu, \tau,\tau^*)$ be a PN space with $\tau$
non-necessarily continuous. Let $(S, F)$ be its associated PSM
space, where $F_{p,q}=\nu_{p-q}$. The strong neighborhoods at $p$
are given by $N_p(t)=\{q\in S: \nu_{p-q}(t)>1-t\}=p+N_\theta(t)$.
In particular, the generalized topology is translation-invariant.
By (N1) we have that this generalized topology is
Fr\'{e}chet-separable (as in the case of PSM spaces). The
countable base of $\theta$-neighborhoods is
$\{N_\theta(\frac{1}{n}): n\in \Na\},$ whose elements are clearly
radial and circled, by axioms (N2) and (N4), respectively.
\end{proof}

Note that the generalized topology induced by a PN space
$(S,\nu, \tau, \tau^*)$ is derivable from the following
generalized uniformity:
$$
\U:=\{V\subset S\times S\mid \exists n\in \Na: V\supseteq \{(p,q)
\mid \nu_{p-q}(\frac{1}{n})\geq 1-\frac{1}{n}\}\},
$$
which is translation-invariant and has a countable filter base of
radial and circled vicinities.

Adapting the methods in \cite{Hoh}, we next show that a converse
result holds for such generalized topologies (or generalized
uniformities).

Let $S$ be a vector space and $(\U_p)_{p\in S}$ be a Fr\'{e}chet-separated,
translation-invariant, generalized topology
of type $V_D$ on $S$. Then, there is a unique translation-invariant,
Hausdorff-separated, generalized uniformity, which is defined as
follows
$$
\U:= \{V\subseteq S\times S\mid \exists U\in \U_\theta:
V\supseteq\{(p,q) : p-q\in U\}\}.
$$
The analogous result of Theorem \ref{Hohle main theorem} for PN
spaces is the following. (Note that there is an extra assumption
on $T$):

\begin{theorem}
\label{adapted-Hohle-theorem} Let $T$ be a $t$-norm such that
$\sup_{0\leq x<1} T(x,x)<1$. Suppose that $T(x,y)\leq xy$,
whenever $x, y<\delta$, for some $\delta>0$. A Fr\'{e}chet-separated,
translation-invariant, generalized topology
$(\U_p)_{p\in S}$ on a real vector space $S$ is derivable from a
Menger PN space $(S,\nu,\tau_T,\tau_{T^*})$, if and only
$\U_{\theta}$ admits a countable base of radial and circled
subsets.
\end{theorem}
\begin{proof}
The direct implication has been shown above. For the converse, let
$\B=\{V_n \mid n\in \Na\}$ be a countable filter base for
$\U_\theta$ consisting on radial and circled
$\theta$-neighborhoods.

Let $N_0\in \Na$ such that $1-\frac{1}{N_0}\geq \sup_{0\leq x<
1}T(x,x)$. We can assume that $\frac{1}{N_0}< \delta$, so that
$T(x,y)\leq xy$, for all $x, y\leq \frac{1}{N_0}$, where $\delta$
is given by hypothesis.

Before defining $\nu$, recall from \cite[Theorem~2]{Hoh} the
distribution functions $F_n$ (used to define the probabilistic
metric $F$):
$$
F_n(x):= \left\{
\begin{array}{ll}
 0 & : x\leq 0 \\
 1-1/(N_0(n+1))  & : 0<x\leq \frac{1}{n+1} \\
 1-1/(2N_0(n+1))  & : \frac{1}{n+1}<x\leq 1 \\
 1-1/(2^{m+1}N_0(n+1))  & : m <x\leq m+1\quad \mbox{for $m\in \Na$}.
\end{array}\right.
$$
By putting ``$\nu_p=F_{p,\theta}$'' in \cite[Theorem~2]{Hoh}) we
define:
$$
\nu_{p}(x):= \left\{
\begin{array}{ll}
F_0  & : p\not\in V_1 \\
F_n  & : p\in V_n\setminus V_{n+1},  \mbox{for $n\in \Na$} \\
\varepsilon_0  & : p\in \cap_n V_n.
\end{array}\right.
$$
We next check that $(S, \nu, \tau_T, \tau_{T^*})$ is a PN space.
Axiom (N1) holds because the generalized topology is
Fr\'{e}chet-separable. (N2) holds because all $V_n$'s are radial. (N3)
holds as in \cite{Hoh}:
$$\tau_T(\nu_p, \nu_q)(x)=\sup_{r+s=x} T(\nu_p(r), \nu_q(s))\leq
1-1/N_0\leq \nu_{p+q}(r+s) =\nu_{p+q}(x).
$$
Finally, for (N4): Let $p\in V_n$ and $\lambda\in [0,1]$. Then,
$\lambda p$ and $(1-\lambda) p$ are also in $V_n$, because $V_n$
is circled. For $x=r+s$, we have to show that $$\nu_p(x)\leq
T^*(\nu_{\lambda p}(r), \nu_{(1-\lambda) p} (s)).$$ Suppose first
that $r$ and $s$ $\geq 1$. Let $a, b, c\in \Na$ such that $a<
r\leq a+1$, $b<s\leq b+1$, and $c<r+s\leq c+1$. Then,
$$
\begin{array}{rcl}
\nu_{\lambda p}(r) &= & 1- 1/(2^{a+1}N_0(n+1)),\\
\nu_{(1-\lambda) p}(s) & = &1- 1/(2^{b+1}N_0(n+1)), \\
\nu_{p}(r+s) & = & 1- 1/(2^{c+1}N_0(n+1)).
\end{array}
$$
By the properties of $T$ it follows that
$$
\begin{array}{rl}
T^*(\nu_{\lambda p}(r), \nu_{(1-\lambda) p} (s))= &
1-T(1-\nu_{\lambda p}(r), 1-\nu_{(1-\lambda) p} (s))\\
=& 1-T(1/(2^{a+1}N_0(n+1)), 1/(2^{b+1}N_0(n+1)))\\
\geq & 1-(1/(2^{a+1}N_0(n+1))) \cdot (1/(2^{b+1}N_0(n+1)))\\
\geq & 1-1/(2^{c+1}N_0(n+1))\\
= & \nu_p(r+s) =\nu_p(x).
\end{array}
$$
In the third line we have used that the arguments of $T$ are
smaller than $1/N_0$, thus we can apply $T(x,y)\leq xy$. Then, we
obtain $\nu_p\leq \tau_{T^*}(\nu_{\lambda p}, \nu_{(1-\lambda
p)})$ as desired. The inequality for the other possible values of
$r$ and $s$, is checked in a similar way. We conclude that $(S,
\nu, \tau_T, \tau_{T^*})$ is a Menger PN space under~$T$.

It only rests to show that the generalized topology induced by
$\nu$ is the same as the one given at the beginning. As in
\cite{Hoh}, we have by construction that
$$
V_n=\left\{p\in S\mid \nu_p\left(\frac{1}{n+1}\right)\geq
1-\frac{1}{N_0(n+1)}\right\}.
$$
Thus, the filter base $\{p\in S\mid \nu_p(\frac{1}{n+1})\geq
1-\frac{1}{n+1}\}$ induced by $\nu$ is equivalent to $\B$, hence
the proof is finished.
\end{proof}
\begin{remark}
{\rm Theorem \ref{maintheorem} also holds if instead of assuming
$T(x,y)\leq xy$ near the origin, one assumes that $T$ is
Archimedean near the origin (i.e.~there is a $\delta>0$ such that
$0<T(x,x)<x$, for all $0<x<\delta$). In that case, the
distribution function $F_n$ can be chosen as:
$$
F_n(x):= \left\{
\begin{array}{ll}
 0 & : x\leq 0 \\
 1-z  & : 0<x\leq \frac{1}{n+1} \\
 1-T(z,z)  & : \frac{1}{n+1}<x\leq 1 \\
 1-T^{m+1}(z,z) & : m <x\leq m+1\quad \mbox{for $m\in \Na$},
\end{array}\right.
$$
where $z=1/(N_0(n+1))$, $T^1(x,y)=T(x,y)$ and recursively
$$T^r(x,y)=T(T^{r-1}(x,y), T^{r-1}(x,y)).$$ }
\end{remark}


\end{document}